\documentclass[a4paper]{amsart}
\usepackage[applemac]{inputenc}
\usepackage{picinpar}
\usepackage[oztex]{graphicx}

\newtheorem{Teorema}{Theorem}

\numberwithin{equation}{section}
\def\sen{\operatorname{\text{sin}}}
\def\Li{\operatorname{\text{Li}}}

\def\Re{\operatorname{\text{Re}}}
\def\Im{\operatorname{\text{Im}}}

\def\R{{\mathbf R}}
\font\cmbsy cmbsy10 at 12 pt

\def\Orden{\mathop{\hbox{\cmbsy O}}\nolimits}   % simbolo de Landau
\begin{document}

% Topmatter
\title[X-Ray of Riemann's Zeta-function]%
   {X-ray of Riemann's Zeta-function}
\author{J. Arias-de-Reyna}
\address{Facultad de Matemáticas\\
    Universidad de Sevilla \\
    P.~O.~box 1160\\
    41080-Seville\\
    Spain}
\email{arias@us.es}
\thanks{Research supported by the Spanish Governement under grant BFM2000-0514.}
\date{September 13, 2002}

\maketitle

\begin{large}

\section{Introduction} \label{S:intro}

This paper is the result of the effort to give the students of the subject {\sl
Analytic Number Theory} an idea of the complexity of the behaviour of the
Riemann zeta-function. I tried to make them {\bf see} with their own
eyes the mystery contained in its apparently simple definition.

There are precedents for the figures we are about to present. In the tables of
Jahnke-Emde \cite{JE} we can find pictures of the zeta-function and some other graphs
where we can see  some of the lines we draw. 
In the dissertation of A. Utzinger \cite{U}, directed by Speiser, 
the lines $\Re\zeta(s)=0$ and $\Im\zeta(s)=0$ are drawn on the 
rectangle $(-9,10)\times(0,29)$.

%I have not been able to consult the
%dissertation of A. Utzinger \cite{U} where, according to Speiser \cite{Spe}, the lines
%where the zeta-function is real are drawn.

Besides, Speiser's paper contains some very interesting ideas. He proves
that the Riemann Hypothesis is equivalent to the fact that the non trivial
zeros of  $\zeta'(s)$ are on the right of the critical line. He
proves this claim using an entirely geometric reasoning that is on the
borderline between the proved and the admissible. Afterwards rigorous proofs of
this statement have been given.
\medskip

Our figures arise from a simple idea. If $f(z)=u(z) + i v(z)$ is a meromorphic
function, then the curves $u=0$ and $v=0$ meet precisely at the zeros and 
poles of the function. That is the reason why we mark the curves where the
function is real or the curves where it is imaginary on the $z$-plane. In order
to distinguish one from the other, we will draw with thick lines the curves
where the function is real and with thin lines the curves where the function is
imaginary.

When I distributed the first figure (the X ray of the zeta function) to my students, I was
surprised at the amount of things one could see in the graphic. I spent a whole
hour commenting on this figure.

Afterwards, I have kept thinking about these graphics. I believe they can be
used to systematize the knowledge that today is scattered. The graphics give
it a coherency which makes it easier to remember.
\goodbreak

\section{The X ray}

\subsection{Remarkable Points} 
Recall that the Riemann zeta-function is defined
for
$\Re (s) >1$ by the series
$$\zeta(s)=1+{1\over 2^s}+{1\over 3^s}+\cdots$$
but it is possible to extend it to the whole plane as a meromorphic function
with a single pole in $s=1$, which is simple.

We already have the X ray on sight. The main thing to take into account
is that the thick lines are formed by those points $s$ where $\zeta
(s)$ is real, and the thin lines by those where $\zeta(s)$ is imaginary.

In the figure we see that, in fact, the lines have a simpler behaviour on the
right of the line $s=1$, that is, on the right of the grey strip which
represents the so called {\bf critical strip}: $0\leq \Re (s) \leq 1$.

In the figure {\bf the real axis}, {\bf the oval} and a thick line that
surrounds the oval stand specially out.

We can see that the real axis cuts lots of thin lines, first the oval
in the pole $s=1$, and in $s=-2$, which is a zero of the function, later, a
line in $s=-4$, another in $s=-6$, \dots, which are the so called {\bf trivial
zeros} of the zeta-function. In the figure we can see these zeros up to $s=-28$,
because the graph represents the rectangle $(-30, 10)\times (-10, 40)$. These
zeros situated in the negative part of the real axis seem to draw the spine of
this X ray.

In the critical strip we see that thick and thin lines meet, that is, we detect
the existence of non trivial zeros. With two decimal digits they are $0.5+ i 14.13$, $0.5 + i 21.02$,
$0.5+i 25.01$, $0.5 + i 30.42$,
$0.5+i 32.93$, $ 0.5+i 37.58$. 
They seem to have a real part equal to ${1 \over 2}$. In fact, it is
possible to prove it.
\medskip

The next remarkable points are the points in which the real axis meets a thick
line. These are points at which the derivative $\zeta'(s)=0$. Since the function
$\zeta(s)$ is real when $s$ is real and since it has zeros in the points $-2n$,
Rolle's theorem from elemental calculus tells us that these zeros have to exist,
at least one between each pair of consecutive even numbers.

For every meromorphic function $f(s)=u(s) + i v(s)$, the lines $u=0$ (and
the lines $v=0$) are smooth curves, except in the points where the derivative
vanishes, in which case two or more curves meet.

An important property that is valid in the general case of a meromorphic
function is the {\bf monotonicity} of the function along the curves.

Let us think of a thick line $v=0$. In this curve, the function $f$
is real. At each point of this curve, the derivative of $v$ in the 
direction of
the curve vanishes. If the derivative of $f$ does not vanish in any point of
the curve, neither will the derivative of $u=f$.

\end{large}

\begin{figure}
\includegraphics[width=\hsize]{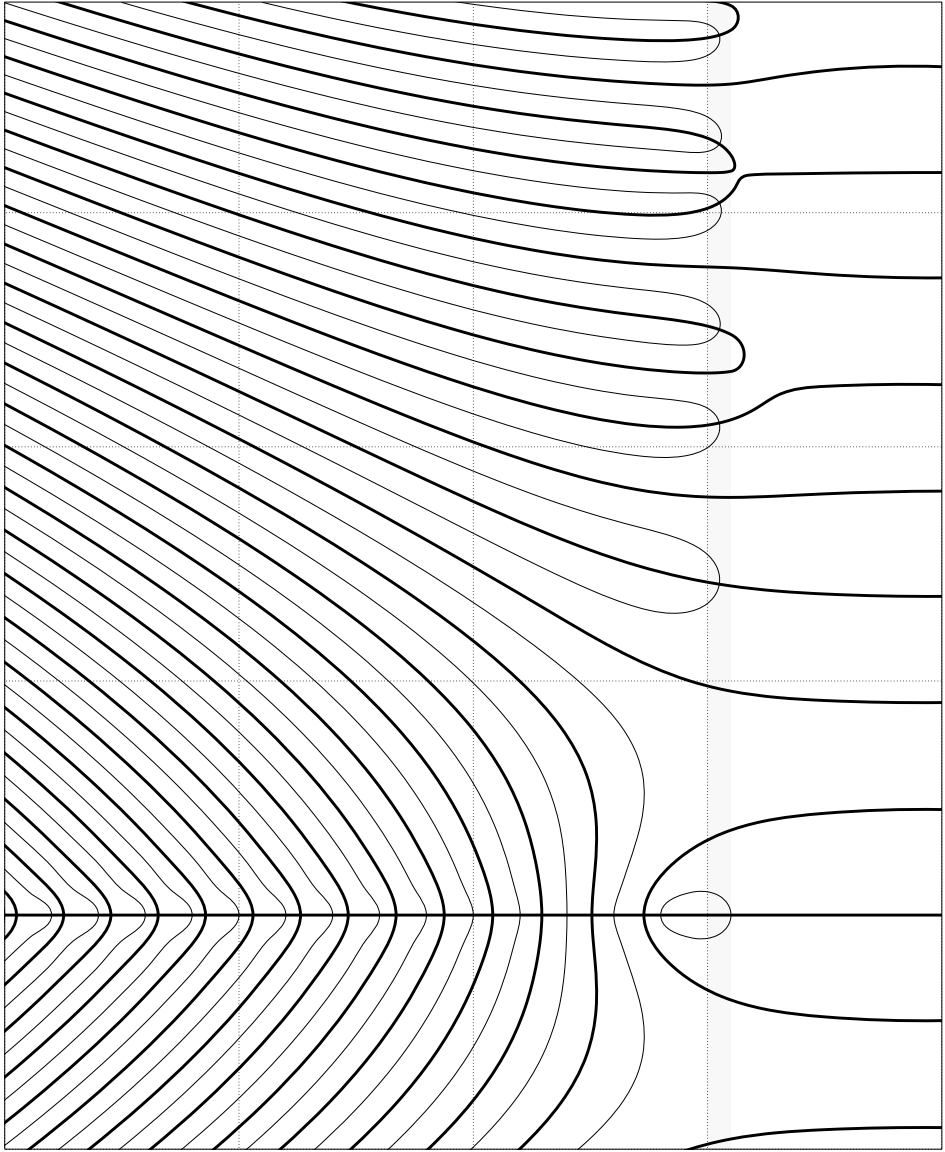}
\caption{X ray of  $\zeta(s)$}
\end{figure}

\begin{large}

So, if on a portion of a curve the derivative of $f$ does not vanish, then the
function varies monotonically as we moves along the curve. In fact, the derivative of
$u$ in the direction of the curve can not vanish on that portion of the curve, so
that it will keep having a constant sign.

When we are dealing with the $u=0$ curves, the monotonous function is $v=f/i$.
\medskip

We would like to know, for each line and each direction, whether the function
increases or decreases. But before tackling this point, we must give each line
a name.

\subsection{The numbering of the lines.} 

 There are three special lines that we will
not number: the oval, the one that surrounds the oval, and the real axis.

The rest of the lines can be numbered with integer numbers. In the first place,
there are thin lines that pass through the zeros in $ -4$, $-6$, $-8$, \dots Line
$-2n$ is the line which cuts the real axis at the point $-2n$.

Between these lines we can find thick lines. For example, there is one line
that goes between the lines numbered $-4$ and $-6$. This one does not go
through the point $-5$, but it cuts the real axis on a point that is between
$-4$ and $-6$. We will say that this one is line $-5$. In the same way lines $-(2n
+ 1)$, $n=2$, $3$, \dots, are defined. These are thick lines that cut the real
axis between $-2n$ and $-2n - 2$.

In this way, we now know which are the lines numbered $n$ for $n=-4$, 
$-5$, $-6$,
\dots. If we consider now how these lines cut the left border of the figure
(line
$x=-30$), we see that the lines already numbered are followed by
other lines, so that it seems natural to call line $-3$ the thick line which
runs parallel to line $-4$, above it, until the latter turns to cut the real
axis while the thick line goes on to the right, crossing the critical strip
at a height near $10$.

So, the even lines are always lines where the function $\zeta(s)$ is purely
imaginary, while the odd lines are thick lines in which $\zeta(s)$ is real. We
must also add that, in many cases, two lines join in a zero. For example, lines
$-2$ and $0$ join in the first non trivial zero of the zeta-function. Lines
$5$ and $7$ also join in the third zero.

This numbering does not include the symmetrical lines below the real axis.

Now we can speak of concrete lines. Do you see how line $11$ turns? It seems it
intended to continue with line $13$ or $15$, but in the end it follows another
path, cutting lines $10$ and $12$.

\subsection{Lines on the right of the critical strip.} 
The
behaviour on the right of the critical strip is governed by the
fact that there are no thin lines. This is easy to understand. For
$\Re (s)=\sigma \to +\infty$ we have
$\zeta(s)\to 1$ uniformly. So, there exists an abscissa $\sigma_0$
so that for $\sigma >\sigma_0$ it holds $\Re \zeta(s)>0$. It
follows that the function does not take purely imaginary values in
this half-plane. It is not difficult to prove that there are no
thin lines passing through the half-plane $\sigma>1.63622\dots$.

\medskip

On the right, the only thick lines that exist are essentially parallel to the
$x$ axis, and they are equally spaced. In order to understand the reason, 
\goodbreak

\end{large}

\begin{figure}
\includegraphics[width=\hsize]{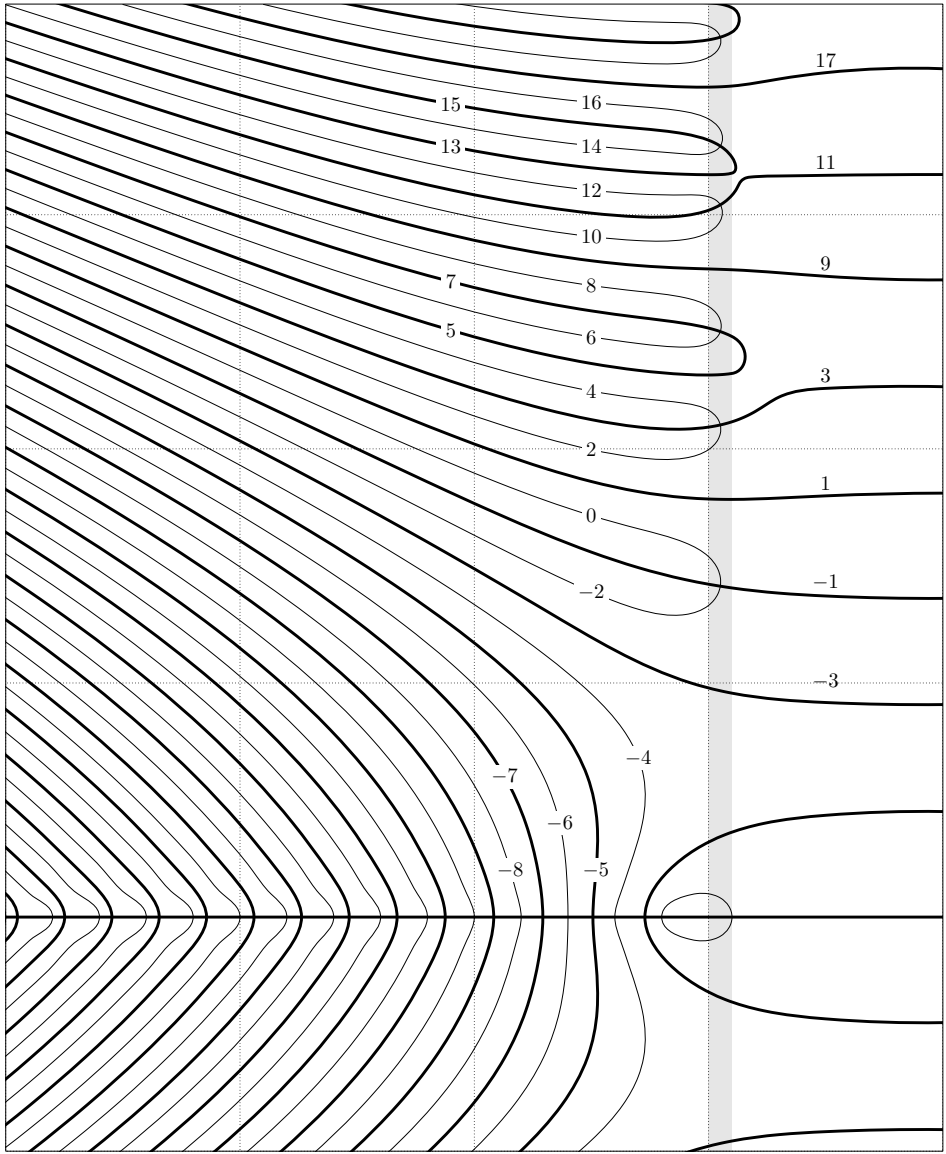}
\caption{Numbering the lines}
\end{figure}

\begin{large}

\noindent 
we start with 
$$\Im\zeta(s)=-\sum_{n=2}^{\infty}{\sen(t\log n)\over n^\sigma}.$$

For big enough values of $\sigma$, the first term of this sum dominates the
rest of them. It is clear that this first term vanishes for $t=n\pi/\log 2$.
So, some parallel lines, separated by a distance which is
approximately equal to $\pi/\log 2 \approx 4.53236014\dots$, exist.

We can see in the figure that these parallel lines, when crossing the critical
strip, alternatively, contain a non trivial zero, or not. It is easy to
explain. The derivative of $\zeta(s)$ along the line which goes at a height
$n\pi/\log 2$, when $\sigma \gg 0$, is given by
$${\partial\over\partial \sigma}\zeta(s)=
{\partial\over\partial \sigma}\Bigl(1+{\cos(t\log2)\over 
2^\sigma}+\cdots\Bigr)\approx -{\cos(t\log 2)\over2^\sigma}\log 
2\approx -{(-1)^n\over 2^\sigma}\log 2.$$
Thus, when the function $\zeta(s)$ runs along this curve from
right to left, (parting from $1$), it is increasing for even values
of $n$. So, for even $n$, $\zeta(s)$ will take on this line the
values in $(1, +\infty)$. On the other hand, for odd $n$ it will
take the values in $(1, -\infty)$, and, in particular, it will
vanish.

Of course we are using the fact that the function is monotonous on the lines,
which depends on the fact that the derivative $\zeta'(s)$ does not vanish on
these lines. This is what happens to the line surrounding the oval, which
ought to contain a zero but does not. So we shall call one of these lines 
{\em parallel}
only when  $\zeta'(s)\ne0$ on it.

\subsection{Orientating ourselves.} 
It holds $\lim_{\sigma\to+\infty}\zeta(\sigma+it)=1$.
 uniformly in $t$. Let us suppose
that we go along line $-1$, starting form the zero on the critical
line, and going to the right. The function will take values in the
interval $(0, 1)$, starting from zero and tending to $1$. The
points we leave on our left, that is, the ones which are a little
above line $-1$, will map to points on the left of the segment
$(0, 1)$, that means, points from the first quadrant.

If, on the contrary, we start from the zero and move to the left
following line $-1$ between lines $-2$ and $0$, the zeta-function
will take negative values, in the interval $(-\infty, 0)$. The
points between lines $-1$ and $0$ will map to points from the
second quadrant.

If we situate ourselves at the zero with our arms outstretched, the right arm
to the right and the left one to the left, we see in front of us a thin line
(line $0$), on which the function will take values $ix$ with $x>0$. Behind us
we have a line (line $-2$) in which it will take values $ix$ with $x<0$. On the
northeast we have points $s$ that map to values $\zeta(s)$ situated on the
first quadrant. On the northwest we have a region which maps into the second
quadrant, etc.

In a line such as line $-3$, which comes from the right and does
not contain any zeros, the zeta-function takes real values greater
than $1$, precisely all the points in the interval $(1, +\infty)$,
because it is known that, on the left, the modulus of the
zeta-function tends to infinity.

\subsection{The regions of the plane.}
 The lines in the graphics divide the plane
into regions. It is easy to realize that the points of one of these regions have
to map by the function $\zeta(s)$ into points of the same quadrant.

\end{large}

\begin{figure} 
\includegraphics[width=\hsize]{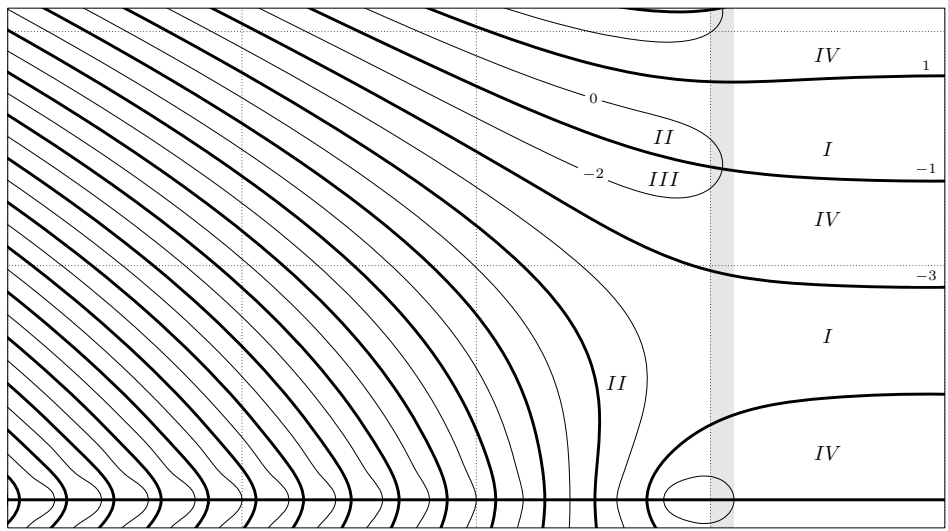}
\caption{Each region transform into a quadrant.}
\end{figure}

\begin{large}

If the points of a region $\Omega$ map to the first quadrant, for example,
and crossing a thin line we reach another region, then the points of this new
region will map to points of the second quadrant; if, on the contrary, we
leave $\Omega$ crossing a thick line, then the new region will map to the
fourth quadrant.

We have various ways to know into which quadrant will a given region 
map. For
example, if we situate ourselves at a non trivial zero and we orientate in the
way described above, then the region in the northeast will map into the first
quadrant, the one situated in the northwest in the second, and so on.

Another way to see it could be to take into account that, if we
walk along a thick line, so that the function $\zeta(s)$, which is
real, increases, then we know that on our left we have a region
that maps into the first quadrant and on our right one that 
maps into the fourth quadrant. For example, let us consider the
thick line number $1$ (this is the line that  comes from the right
and goes between the first and the second non trivial zeros of the
zeta-function). and let us walk through it to the left. That is,
we walk through this line coming from $+\infty$. The values of the
zeta-function on the points of this line are greater than $1$,
and, as we walk through it, the values on its points are growing.
Then, the region on our right will map into the first quadrant.

\vfil

As an exercise, we frame this question: Into what do the points of the 
oval map?
\goodbreak

\end{large}

\begin{figure}
\includegraphics[width=\hsize]{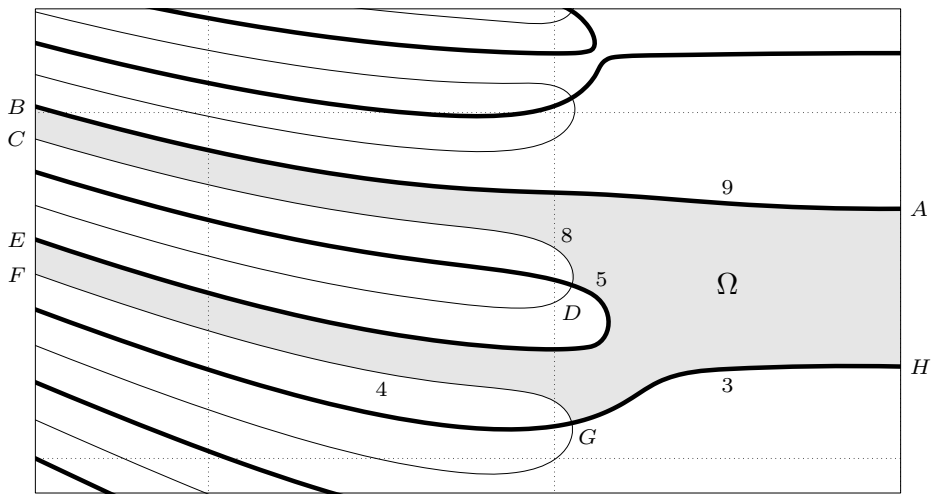}
\caption{The region $\Omega$.}
\end{figure}

\begin{large}

\subsection{The equation $\zeta(s)=a$.}
Let us consider an equation of the form
$\zeta(s)=a$. Where are its solutions? How are they distributed?

If $a$ is situated in a determined quadrant, say the first, then there will
not be any solution in the region bordered by lines $-2$ and $-1$. In fact, the
points of this region map to points from the third quadrant. Thus, the
solutions have to be situated in the regions we know to map into the first
quadrant. We can say how many solutions are there  in each of these regions.

For example, let us consider the region $\Omega$, bordered by lines $9$, $8$,
$5$, $4$ and part of line $3$.

As $s$ runs through the border of this region, ABCDEFGHA, the image $\zeta(s)$
runs through a closed path, according to the following table
$$\vbox{\halign{#\hfil\quad &\hfil#\hfil\quad& \hfil#\hfil\quad& \hfil#\hfil\quad& 
\hfil#\hfil\quad& \hfil#\hfil\cr
$s$&$A\to B$& $C\to D$& $D\to E$& $F\to G$&$G\to H$\cr
\noalign{\medskip\hrule\medskip}
$\zeta(s)$&$1\to+\infty$&$+i\infty\to0$&$0\to+\infty$&$+i\infty\to0$&$0\to1$\cr}}$$
so that $\zeta(s)$ goes round the first quadrant twice. Thus, if $a$ is a point
from the first quadrant, the equation $\zeta(s)=a$ will have exactly two
solutions when $s\in \Omega$, according to the argument principle.

It is clear that, under the former conditions, the equation 
$\zeta'(s)=0$ must
have a solution in the region $\Omega$. There are other solutions of
$\zeta'(s)= 0$, one for each of the lines $-3$, $-5$, $-7$, \dots precisely in the
points at which these lines meet the real axis.

In the X ray figure we thus see $14$ solutions of $\zeta'(s)=0$, situated on
the real axis, and, besides, the existence of two more can be inferred, one
in the region called $\Omega$ and the other in the region situated between
lines $11$ and $17$.

\subsection{The construction of the graphics.}
The first one who
calculated zeros of the function $\zeta(s)$ was Riemann himself.
Gram \cite{G} calculated the first ten zeros and proved that they were
the only ones satisfying $0< \Im (\rho )<50$. The following step
in this direction were the works of Backlund \cite{Ba}, who managed to
prove that the zeros $\alpha + i\beta$, with $\alpha>0$ and
$0<\beta<200$ where exactly $79$, every one of them with real part
$\alpha=1/2$. Later, Hutchinson \cite{H} managed to extend these
calculations up to $300$. This task has been followed by many
mathematicians such us Titchmarsh, Comrie, Turing, Lehmer, Rosser,
and so on.

We have used the same techniques they used to localize certain values of
$\sigma$, the solutions of $\Re\zeta(\sigma + it)= 0$ or $\Im\zeta(\sigma +
it)=0$. Afterwards, we have written other programs, based on these, 
ones which
situate the points belonging to the same line in order.

These graphics owe much to the program Metafont, created by Knuth to design the
fonts used by \TeX, and to its modification MetaPost, by J.~D.~Hobby, which
allows one to obtain PostScript graphics.

For $|\sigma - 1/2|>5/2$ I have proved that the curves follow the
path I have drawn here, but for the area near the critical strip I just
followed the method until I was convinced that the curves follow the path I
have drawn. This has required, in some cases, the calculation of numerous
points on each line. So they can not be considered proved.

\subsection{ The existence of non trivial zeros.}
One of the problems which that arose
in drawing these graphics was to know the number of a line. Applying the
Stirling series for $\log\Gamma(s)$ and the functional equation of $\zeta(s)$ we
can prove the following:
\begin{Teorema}\label{T:A} 
The number of a line passing through the point $-1 + it$ (with
$t>5)$ is the integer number that is nearest
\begin{equation}
{2 t\over\pi}\log{ t\over 2\pi}-{2 t\over \pi}+{1\over2}.\label{E:Uno}\end{equation}
\end{Teorema}
Logically, the role of line $\sigma=-1$ is not important as long as we are not
near the critical strip.

If we watch all the lines coming to the critical strip from the left and we
know that the thin lines can not surpass a certain point, we see they must come
back and they can not do it unless they cut the thick lines that accompany
them. This little rigorous reasoning can be turned into a proof, by means of
the argument principle, of the following theorem.

\vfil\eject

\begin{Teorema}
    Let us suppose that $-1 + iT$ is on one of the
parallel lines which do not contain any zero. The number of non
trivial zeros $\rho=\beta + i\gamma$ below this line and
satisfying $\beta>0$ is equal to the integer number nearest to
\begin{equation}
{T\over2\pi}\log{T\over2\pi}-{T\over2\pi}+{7\over 8}.
\label{E:Dos}
\end{equation}
\end{Teorema}

\begin{minipage}{4cm}
\hspace{-0.8cm}
\includegraphics[width=4cm]{figure05.ps}
\ \vspace{0.5cm}  \  
\end{minipage}
\hfill
\begin{minipage}{7.6cm}
\begin{proof}[Proof]
    In order to prove it, let us consider the region
bordered by the lines $\sigma=-1$, $\sigma=2$, line number $-3$ and
the line which is referred to in the statement. In the figure we
have represented the case when this line is line number $17$. The
number of zeros we are looking for is equal to the variation of
the argument of $\zeta(s)$ along this curve. But, along the
segment $AB$, there is no variation, because all the segment maps
into points of the first and the fourth quadrant, starting and
finishing on the real axis. In the portions of the curves $CB$ and
$DA$ there is also no variation because on them the function $\zeta(s)$ 
 is real positive. Thus, the variation of the argument is exactly the
one which takes place in the segment $CD$. This segment cuts
several lines. As it goes from one to the other, the argument of
$\zeta(s)$ varies precisely in $\pi/2$. This way, the variation of
the argument will be $(N + 3)\pi/2$ if the highest line is line
number $N$. By the previous theorem $N$ is the integer number
nearest to \ref{E:Uno}. It follows easily that the number of zeros is
equal to $(N + 3)/4$, and thus it is the integer number nearest to
\ref{E:Dos}.
\end{proof}

If the line we are dealing with is one of the parallel line which do contain
zeros, an analogous reasoning proves that the number of zeros below this line,
counting the one which is on line $N$ on, is $(N + 5)/4$, and so it is the
integer number which is nearest to
\begin{equation}
{T\over2\pi}\log{ T\over 2\pi}-{T\over 2\pi}+{11\over8}.
\end{equation}
\vfill
\end{minipage}
\noindent 
So, the parallel lines which do not contain zeros have a number
$N\equiv 1$ 
and the ones which do contain zeros, on the contrary, satisfy
$N\equiv 3\pmod{4}$.
\goodbreak
%\vfill

\section{Techniques to calculate the function $\zeta(s)$.}

In order to construct the graphics the possibility of calculating the function is
essential. There are two basic ways to calculate the function $\zeta(s)$.

\subsection{Euler-MacLaurin Formula.}
It is the following 
$$
\zeta(s)=\sum_{n=1}^{N-1}{1\over n^s}+{1\over2}{1\over N^s}+{N^{1-s}\over
s-1}
+
\sum_{k=1}^M T_k+R(N,M),$$
where
$$T_k={B_{2k}\over(2k)!}N^{1-s-2k}\prod_{j=0}^{2k-2}(s+j),$$
$B_n$ are the Bernoulli numbers, and there are good known bounds for the
error term $R(N, M)$.

Choosing $N$ and $M$ conviniently as a function of $s$ allows one to calculate
$\zeta(s)$ for any $s$ and with arbitrary precission.

\subsection{Riemann-Siegel Formula.}
Before explaining what does it consist of, we
should define Hardy's function (called so though it was known by Riemann)
$$Z(t)=e^{i\theta(t)}\zeta\Bigl({1\over2}+i\,t\Bigr)$$
with
$$\theta(t)=\Im\Bigl(\log\Gamma\Bigl({1\over4}+ i\, {t\over2}\Bigr)\Bigr)-{t\over 2}\log
\pi.$$

The functional equation implies that the function $Z(t)$ is real
for real values of $t$. This allows one to locate zeros
$\rho=\beta + i\gamma$ with $\beta=1/2$, because a change in the
sign of $Z(t)$ implies the existence of a zero with abscissa
exactly $1/2$.

The Riemann-Siegel formula allows one to calculate the function for any point
$s$, but with a limited precission depending on who $s$ is. We will only write
it for a point of the form $s=1/2 + it$. It is the following:
$$Z(t)=2\sum_{n=1}^m{\cos(\theta(t)-t\log n)\over \sqrt{n}}+g(t)+R,$$
where

\noindent
\begin{minipage}{3cm}\ \hskip1cm
\parbox{3cm}{\begin{eqnarray}
m&=&\lfloor\sqrt{t/2\pi}\rfloor\nonumber\\
g(t)&=&(-1)^{m-1}\Bigl({t\over2\pi}\Bigr)^{-1/4}h(\xi)\nonumber\\
h(\xi)&=&(\sec 2\pi \xi)\cos 2\pi\phi\nonumber
\end{eqnarray}}
\end{minipage}
\hfill
\begin{minipage}{3cm}
\parbox{3cm}{\begin{eqnarray}
\xi&=&\Bigl({t\over 2\pi}\Bigr)^{1/2}-m\nonumber\\
\phi&=&\xi-\xi^2+1/16\nonumber\end{eqnarray}
}\hskip1cm\ \end{minipage}
\goodbreak
\bigskip

The error is $\Orden(t^{-3/4})$. The Riemann-Siegel formula requires, in
order to achieve a fixed precission, the calculation of approximately
$\sqrt{t}$ terms of the first sum, while the Euler-MacLaurin formula requires
a number of terms of the order of $t$. When $t\approx 10^9$, this difference
turns out to be very significant.

This formula has an interesting story. In a letter directed to Weierstrass,
Riemann claims that:
\bigskip

\begin{minipage}{10cm}
the two theorems I have just stated here: {\sl 
\medskip
\noindent{\em 
that between $0$ and $T$ there
exist approximately ${T\over2\pi}\log{T\over2\pi}-{T\over 2\pi}$
real roots of the equation $\xi(\alpha)=0$.} 

\medskip
\noindent{\em 
that the series $\sum_{\alpha}\bigl(\Li(x^{1/2+i\alpha})+\Li(x^{1/2-i\alpha})\bigr)$,
when its terms are ordered according to the incresing
growth of $\alpha$, tends to the same limit than the expression 
$${1\over2\pi i\log x}\int\limits_{a-bi}^{a+bi}{d{1\over
s}\log{\xi\bigl((s-{1\over2})i\bigr)\over\xi(0)}\over ds}x^s\,ds$$
when $b$ grows beyond boundaries.}

\medskip
\noindent are consecuences of a new expansion $\xi$ which I have not been able to simplify
enough to comunicate it.}
\end{minipage}

\bigskip
Riemann considered proved\footnote{It seems that the 
Riemann's capability to prove his claims has been doubted
 more than once without proper support. To get more information on this point
you can consult my paper\cite{A}.} all his other claims contained in \cite{R}.

Because of this, it was clear that, among Riemann's papers should be the
expansion he was talking about. Some $70$ years after his death, 
C.~L.~Siegel
\cite{Si}  managed to figure out these papers. It would be worth for the reader to
have a look at the photocopy of the sheets of Riemann's manuscript which
contain the famous expansion, which can be found in Edward's book \cite{E}, to
realize the difficulty of Siegel's task.

\section{Graphics of the zeta-function in the critical line.}

In the following pages we present some graphics corresponding to the values of
$t$ between $0$ and $560$. We have just drawn the strip $-1\le \sigma \le 2$.
The part which is not represented has nothing new and our imagination can
supply it without any trouble. On the left we have pointed the values of $t$,
each time it increases in twenty unities, and we have written the number of
some lines (the parallel ones which do not contain any zeros). We leave till
later the explanations about some points we have marked in the graphics with
little circles, situated on the critical line.

In the usual proof about the number of non trivial zeros up till a
given height $T$ it is proved that a horizontal segment cutting
the critical strip between $\sigma=-1$ and $\sigma=2$, at a height
of, say, $t$, can only cut a number of lines of the order $\log
t$. This, and the fact that the lines coming from the right can
not cut themselves, because that would generate a region bordered
by points where the function $\zeta(s)$ is real and bounded,
allows Speiser to infer that the lines coming from the right have to
cross the critical strip and move away to the right towards the
infinity.

In fact, if one of these lines went up towards the infinity through the critical
strip it would force the others to go above it, also through the critical line,
to infinity. As there are a number of them of the order of $t$, we would reach
a contradiction with the results already proved about the horizontal segment in
the critical strip.

As the number of lines on the left has an order of $t\log t$, we see that it is
necessary for some of them to go back to the right, joining others. Speiser
calls these figures, formed by two lines, when the function is real, 
{\sl sheets}.

\medskip

The graphics give rise to some integer sequences, the most obvious is the
sequence formed by the number of lines which escape to the right. That is,
$$ -3, -1, 1, 3, 9, 11, 17, 23, 29, 35, 41, 47, 53, 59, 69, 75, 81,
91, 97, 103, 113,$$
$$123, 129, 135, 145, 155, 161, 171, 181, 187, 197, 
207, 217, 223, 237, 247, 253$$
$$263, 273, 283, 293, 307, 313, 323, 329, 343, 353, 359, 373, 383, 
393, 403, 417$$
$$423, 437, 451, 457, 467, 481, 491, 501, 511, 525, 535, 545, 559, 
569, 579,\dots$$
Speiser suggests that it could be connected with the distribution of the prime
numbers. I could  not see any conection.

Looking at these graphics, a general scheme seems to appear. Between two
parallel lines which do not contain zeros there are an even number of thin
lines, joining by pairs, a thick line coming parallel from the right cuts one
of these loops and the rest are inserted with loops of thick lines. Later we
will see how these simple ideas about the function break.

The thin lines almost do not cross over the critical line, while the thick
lines sometimes reach $\Re(s)=2$. In fact, we see that at a high height they do
cross over this line, this happens for the first time in the sheet formed by
lines $789$ and $791$. But what is really surprising happens in line $1085$.
This forms a sheet with line $1091$. So this sheet surrounds completely the one
formed by lines $1087$ and $1089$. For higher values of $t$ this becomes
quite frequent.
\goodbreak

\end{large}
\begin{minipage}{2.9cm} 
\includegraphics[height=\vsize]{figure06.ps}
\end{minipage}
\begin{minipage}{2.9cm}
\includegraphics[height=\vsize]{figure07.ps}
\end{minipage}
\begin{minipage}{2.9cm}
\includegraphics[height=\vsize]{figure08.ps}
\end{minipage}
\begin{minipage}{2.9cm}
\includegraphics[height=\vsize]{figure09.ps}
\end{minipage}
\vfil\eject

\begin{minipage}{2.9cm}
\includegraphics[height=\vsize]{figure10.ps}
\end{minipage}
\begin{minipage}{2.9cm}
\includegraphics[height=\vsize]{figure11.ps}
\end{minipage}
\begin{minipage}{2.9cm}
\includegraphics[height=\vsize]{figure12.ps}
\end{minipage}
\begin{minipage}{2.9cm} 
\includegraphics[height=\vsize]{figure13.ps}
\end{minipage}
\vfil\eject

\begin{minipage}{2.9cm}
\includegraphics[height=\vsize]{figure14.ps}
\end{minipage}
\begin{minipage}{2.9cm}
\includegraphics[height=\vsize]{figure15.ps}
\end{minipage}
\begin{minipage}{2.9cm}
\includegraphics[height=\vsize]{figure16.ps}
\end{minipage}
\begin{minipage}{2.9cm} 
\includegraphics[height=\vsize]{figure17.ps}
\end{minipage}
\vfil\eject

\begin{minipage}{61.6pt}
\includegraphics[height=550.4pt]{figure18.ps}
\end{minipage}\hfil
\begin{minipage}{60.8pt}
\includegraphics[height=549.6pt]{figure19.ps}
\end{minipage}\hfil
\begin{minipage}{171.64125pt}
\begin{minipage}{171.64125pt}
\includegraphics[height=170.pt]{figure20.ps}
\medskip

\begin{minipage}{171.64125pt}
    Sometimes, as in line $1187$, can happen that the path followed by the line is
not clear. In such cases we draw the confusing area, for example, in this case
we can see here a square which represents the rectangle $(0.2, 0.8)\times (540,
540.8)$, marking the points where the function $\zeta(s)$ is real. This way it
is clear that line $1187$ turns down when it gets to the critical line, and the
pair of lines $1189$ and $1191$ form a sheet.
\bigskip

The monotonicity of the function on the curves is valid except in the poles,
provided we define suitably the prolongation of the curve in the points where
the derivative vanishes. In fact, in one of these points several curves meet,
and there are as many on which the function increases as on which the function
decreases, thus it always is possible, when we reach a zero of order $n$ of the
derivative, to go back with an angle of $2\pi/n$, so that the 
monotonicity is
maintained all through the curve.
\bigskip

Another particularity of the curves we would like to emphasize is that,
except possibly the ones which go through the poles, the curves have to
leave every compact. Thus they come from and go to the infinity.
\end{minipage}
\end{minipage}
\end{minipage}

\section{First Theorem of Speiser}

\begin{Teorema}[Speiser]
The Riemann Hypothesis is equivalent to the fact that
all the sheets meet the critical line.
\end{Teorema}

\begin{proof}[Proof]
    Let us recall that the sheets are the thick lines in which $\zeta(s)$ is
real and which, coming form the left, go back to the left joining another 
thick line.

In the first place, on each sheet, the zeta-function takes all the
real values. In fact, the fuction is real and monotonous and it
tends to infinity as $\sigma \to -\infty$. Thus in a sheet
there is always one (and only one) zero of the function.
\medskip

If some sheet would not touch the critical line, the corresponding zero would
be an exception to the Riemann Hypothesis.
\bigskip

The proof of the other implication is a little more intricate. In the first
place, we must notice that if a real line crosses the critical line through a
point $\alpha$, then it holds that $\zeta(s)=0$ or else the modulus
$|\zeta(s)|$ decreases as the curve goes through this point from left
to right. To see this, we notice that, since the function of Hardy $Z(t)$ is
real for real values of $t$, it follows that $\theta(t) + \arg\zeta(1/2 + it)=
\hbox{cte}$, except  at one zero of the zeta-funcion.
Since $\theta(t)$ is an increasing function it follows that  
$\partial_{t}\arg\zeta(1/2+it)<0$ unless $1/2+it$ is a zero of the 
zeta-function. 
Let us consider then the
analytic function
$$\log\zeta(1/2+it)=\log|\zeta(1/2+it)|+i\arg\zeta(1/2+it).$$
From the Cauchy-Riemann equations it now follows that the derivative of
$|\zeta(s)|$ with respect to $\sigma$ in the point $\sigma=1/2$ must be
negative.
\medskip

If a line (of a sheet) surpasses the critical line from left to right it must
go through it again in order to come back to the left. If the zero which is on
the sheet is not one of the points where the line cuts the critical line, we
can apply the results obtained above to both of them. It follows that in one of
the points it must be satisfied that $\zeta(s)>0$, and in the other,
$\zeta(s')<0$, so that the absolute value can be decreasing in both. Since the
function is monotonous along the sheet, it follows that the zero is situated on
the right of the critical line.

Thus there are two kinds of sheets cutting the critical line. In some of them,
the zero is located in one of the two cuts between the critical line and 
the sheet. In
others, the zero is situated on the right of the critical line.

Consecuently, if the Riemann Hypothesis is false, there would be a zero
situated on the left of the critical line, which can be situated only on a
sheet which does not cut the critical line. If the Riemann Hypothesis is true,
then all the sheets are of the kind which have a zero exactly on the line and
so they cross it (or are tangent at it).
 
This way, we see that if the Riemann Hypothesis is true all the sheets would
meet the critical line.
\end{proof}

In the figures above we can see that, in fact, all the sheets do contain a zero
and another point on which it cuts the critical line. We have marked the latter
with a\break

\noindent  little circle.
These points are the so called Gram points, and we will
see  some interesting results about these points below.
\vfil
\noindent
\begin{minipage}{356pt}
\begin{minipage}{52.185pt} 
\includegraphics[height=505.68pt]{figure21.ps}
\end{minipage}\hfil
\begin{minipage}{52.185pt}
\includegraphics[height=504.945pt]{figure22.ps}
\end{minipage}\hfil
\begin{minipage}{55.86pt}
\includegraphics[height=505.68pt]{figure23.ps}
\end{minipage}\hfil
\begin{minipage}{55.86pt}
\includegraphics[height=505.68pt]{figure24.ps}
\end{minipage}
\centerline{\  }
\centerline{{\scshape Figure} 12. Zeros of  $\zeta(s)$}
\end{minipage}
\eject

\begin{large}
    
\section{Separating the zeros. Gram points.}
In the last page we marked the zeros of the function $\zeta(s)$. We see that
there is some randomness in their distribution.

To try to understand their distribution the Gram points are brought in. The Gram
points, which we have marked with a little circle in our figures, are points
on the critical line where the function $\zeta(s)$ is real and does not
vanish. In the graphics of the following page we have marked Gram points. It
is remarkable how regularly are these points distributed.
\medskip

The Gram point $g_{n}$ is defined as the solution of the following equation:
$$\theta(g_n)=n\pi$$

Since Hardy's function is real,
$$\zeta(1/2+i \,g_n)=e^{-i\theta(g_n)}Z(g_n)=(-1)^n Z(g_n)$$ 
is real. Thus these points are situated on thick lines.

Looking at the preceding figures, we verify that, in most of the cases
represented there, it holds:
$$\zeta(1/2+i \,g_n)>0.$$  
We will see later that this inequality has exceptions for greater values of $t$.

In the Riemann-Siegel formula, the first term, which is the most important, has
a value of $(-1)^n$ for $t=g_n$, which partly explains the tendency of $Z(g_n)$
to have the sign $(-1)^n$.
\medskip

An important consequence of the preceding is that $Z(g_n)$ and $Z(g_{n + 1})$
have oposed signs, so that the function $\zeta(1/2 + it)$ will vanish in the
interval $(g_n, g_{n + 1})$. In Figure $14$ we can see this
graphically.

When one intends to calculate the real zeros of a polynomial, a first task is
to separate the zeros, that is, to find an increasing sequence of values for
$t$ in which the polynomial takes alternatively values of different signs. An
analogous aim is achieved with Gram points.

Gram noticed that these points seemed to separate the zeros of the
zeta-function and claimed that this would be true for not too high values of
$t$. Hutchinson named the fact that in each interval $(g_n, g_{n + 1} )$ would
be a zero of the zeta-function Gram's law.

Titchmarsh uses this idea to prove that there is an infinity  of zeros
on the critical line, proving that the mean value of $Z(g_{2n})$ is positive
and that of $Z(g_{2n + 1})$ is negative.

\end{large}
\begin{minipage}{320pt}
\begin{minipage}{56.8pt} 
\includegraphics[height=550.4pt]{figure25.ps}
\end{minipage}\hfil
\begin{minipage}{56.8pt}
\includegraphics[height=549.6pt]{figure26.ps}
\end{minipage}\hfil
\begin{minipage}{60.8pt}
\includegraphics[height=550.4pt]{figure27.ps}
\end{minipage}\hfil
\begin{minipage}{60.8pt}
\includegraphics[height=550.4pt]{figure28.ps}
\end{minipage}
\centerline{\ }
\centerline{{\scshape Figure} 13. Gram Points.}
\end{minipage}
\begin{large}

\end{large}
\begin{minipage}{320pt}
\begin{minipage}{56.8pt} 
\includegraphics[height=550.4pt]{figure29.ps}
\end{minipage}\hfil
\begin{minipage}{56.8pt}
\includegraphics[height=549.6pt]{figure30.ps}
\end{minipage}\hfil
\begin{minipage}{60.8pt}
\includegraphics[height=550.4pt]{figure31.ps}
\end{minipage}\hfil
\begin{minipage}{60.8pt}
\includegraphics[height=550.4pt]{figure32.ps}
\end{minipage}
\centerline{\ }
\centerline{{\scshape Figure} 14. Gram's points separate the zeros.}
\end{minipage}

\begin{large}

\section{Second Theorem of  Speiser.}

We present here Speiser's proof, and, as we have already said, his methods are
between the proved and the acceptable. Everybody quotes him but nobody
reproduces his theorems. His methods do not seem convincing to me, either, 
though I think his proof is essentially sound. We present it more as a challenge: to turn
it into a proof, filling its gaps. In any case, a flawless proof of a stronger
result can be found in Levinson and Montgomery \cite{LM}.

\begin{Teorema}[Speiser]
    The Riemann Hypothesis is equivalent to the fact
that the non trivial zeros of the derivative $\zeta'(s)$ have a real part $\ge
1/2$, that is, that they are on the right of the critical line.
\end{Teorema}

\begin{proof}[Proof]
    Let us assume that there is a zero $a$ of $\zeta(s)$ on the left
of the critical line. Let us consider the lines of constant argument
$\arg\zeta(s)=\hbox{cte}$ which come from the point $a$. In all of them the modulus
$|\zeta(s)|$ is increasing. Thus, these lines can not cross the
critical line, because at crossing it from left to right, the absolute value of
$\zeta(s)$ ought to decrease.

They also, can not be tangent at the critical line, because the tangency point
would be a point in the critical line where $\arg\zeta(1/2 + it)$ would be
stationary, and this is possible only if it is a zero. But, in the line,
$|\zeta(s)|$ increases, starting from zero, so it can not go through a
zero of the function.

Thus all these lines come back to the left. Some of them go back leaving the
point $a$ below, others leave it above. The line which separates both kinds of
lines must reach a zero of the derivative, which will allow it to go back.
This would be a zero of the derivative on the left of the critical line.

Consequently, if the Riemann Hypothesis is false, we see that there must exist
a zero of the derivative on the left of the critical line.
\medskip

Now let us suppose  that $\zeta'(a)=0$, where $\Re(a)<1/2$. We have to find a
zero of the function on the left of the critical line. We can assume that
$\zeta(a)\not=0$, because if this were the case we would have already finished.
Since the derivative vanishes, there exist two opposite lines, of constant
argument and along which $|\zeta(s)|$ decreases. We follow these two
lines, and we must reach a zero, because $|\zeta(s)|$ decreases. If it
is on the left of the critical line, we have finished, while, in other case, it
is clear that we will reach the critical line.

Our two paths, till they reach the critical line, and the segment from the
critical line they determine, enclose a region $\Omega$. From the point $a$,
two opposite paths also set off, along which $|\zeta(s)|$ increases,
and the argument of $\zeta(s)$ is constant. One of them enters our region
$\Omega$, (because in the point $a$, the borderline of the region has a
well-defined tangent). 
\vfil\eject

%\end{large}
\begin{minipage}{320pt}
\begin{minipage}{97.6pt} 
\includegraphics[height=550.4pt]{figure33.ps}
\end{minipage}\hfil
\begin{minipage}{97.6pt}
\includegraphics[height=549.6pt]{figure34.ps}
\end{minipage}\hfil
\begin{minipage}{101.6pt}
\includegraphics[height=550.4pt]{figure35.ps}
\end{minipage}
\centerline{\ }
\centerline{{\scshape Figure} 15. Zeros $\zeta'(s)$.}
\end{minipage}

%\begin{large}

\noindent
We will follow this path. Since $\vert\zeta(s)\vert$
increases to infinity along this curve, it must leave the region $\Omega$, but
it can not do it across the curves which we used to define it since along them
$|\zeta(s)|<|\zeta(a)|$, and we are considering a curve where
the values of $\vert\zeta(s)\vert$ are greater than $|\zeta(a)|$. It
also can not leave $\Omega$ across the segment from the critical line, because
to cross it from left to right it ought to do it with $|\zeta(s)|$
decreasing.

Thus supposing that there was no zero of $\zeta(s)$ on the left of the critical
line has lead us to a contradiction.
\end{proof}

Speiser's theorem makes the zeros of $\zeta'(s)$ far more interesting. R.~Spira
who has given a complete proof of half of Speiser Theorem
 \cite{Sp} has calculated the first ones (those which have an abscissa less than
$100$), which are represented in the figure.

We see that the real curves (thick ones) seem to be attracted by the zeros,
and each sheet seem to have one zero associated with it, which would justify
their crossing the critical line to approach their corresponding zero. In this
way we  have an insight about the place where the zeros which are beyond
line $113$ (that we have not draw in the figure)  are situated.

We see that a zero of the derivative does explain the behaviour of line $11$.
If the derivative vanished at a point in which the function is $\zeta(s)$ is
real, at this point two thick lines would meet perpendiculary. What happens
here is that the function is almost real in the zero and the curves resemble
the meeting we have described. This is also what happens in line $1187$. If we
see the graphics with the dots we had to do in order to see the path this line
follows, we can verify that the derivative at this point has a zero with an abscissa slightly
greater than $1/2$.
\medskip

Because of Speiser's reasoning we know that a line which comes parallel from
the right at a height $T$ does not rise or fall of level until it crosses the
critical strip at a height greater than $\Orden(\log T)$. Thus, the number of
parallel lines below it is $T\log 2/\pi + \Orden(\log T)$. This lines,
alternatively, contain a zero of $\zeta(s)$ which is not associated with a
sheet, or do not contain it. Thus the zeros which are not associated with a
sheet up till a height $T$ is approximately equal to $T\log 2/2\pi$. According
to this, the number of sheets below a height $T$ is
$${T\over2\pi}\log{T\over2\pi}-{T\over2\pi}+{7\over8}-\Bigl({T\log2\over2\pi}-{1\over2}\Bigr)=
{T\over2\pi}\log{T\over4\pi}-{T\over2\pi}+{11\over8},$$
with an error of the order of $\log T$.

If the remarks we have made concerning the zeros of the derivative $\zeta'(s)$
are true this is the number of zeros of the derivative up till a height $T$.
Spira conjectured this claim and Berndt \cite{Be} has proved it.

\section{Looking higher. A Counterexample to  Gram's Law}

We will now see that most of the regularities in the behaviour of the function
break at a great enough height.

\bigskip

\begin{minipage}{6cm}
\hspace{-0.8cm}
\vspace{-0.3cm}                                 %%   MetaPost
\includegraphics[width=5.1cm]{figure36.ps}   %%   RADIOGRAFIA::04-zeta240/320::Gram.mp 
\ \vspace{-0.cm}  \                             %%
\end{minipage}
\hfil
\hspace{-1.4cm}
\begin{minipage}{6.6cm}
    \vspace{-0.3cm}
    Gram's law claims that between every two consecutive Gram points 
    there exists a zero
of the function $\zeta(s)$. The remark is Gram's \cite{G}, but it was Hutchinson who
named it {\sl Gram's law}, although it was him who found the first
counterexample, which can be seen in the figure. Gram \cite{G} only claimed that
this would happen for the first values of $n$.

The interval $(g_{125}, g_{126})$ does not contain zeros of the zeta-function.
On the contrary, the next interval contains two zeros thus reestablishing the
total count.

Later Lehmer \cite{Le} finds out that the exceptions grow more frequent as $n$
increases. He also notices that, in general, these exceptions consist of a Gram
interval in which there are no zeros, next to another which has two.
\medskip

Possibly the only valid rule was the one formulated by Speiser: The number of
thick lines crossing the line $\sigma= -1$ below a height $T$ is
$${T\over\pi}\log{T\over2\pi}-{T\over\pi}+{1\over4}$$
and the number of Gram points below this height is only half this number. A
thick line which crosses the critical line has to do it through a Gram point or
a zero. The parallel lines, alternatively, contain a Gram point or a zero.
Speiser believes that each sheet uses a Gram point and a zero to enter and exit
the area on the right of the critical line. 
\end{minipage}
\vfil\eject

\noindent
This is true, there exists a
bijectiv map between the zeros and Gram points, but it is the one established
by the fact that they are on the same sheet. It may be convenient to call also
sheet to a parallel line which does not contain zeros and the parallel line
inmediately above it.

This way, another sequence of natural numbers associated with the graphics of
the zeta-function arises. In fact, Gram points $g_{-1}$, $g_0$, 
$g_1$, \dots are
associated with the zeros numbered
$$1, 2, 3, 4, 5, 7, 6, 8, 10, 9, 11, 13, 12, 14, 16, 15, 17, 18, 
20,19,21, 23, 24, 22, 26,$$
$$25, 27, 28, 30, 31, 29, 32, 34, 33, 35, 
36, \dots$$
It is a permutation $\sigma$ of the natural numbers, so that $|\sigma(n) -
n|\le C\log n$. But, actually, it is well defined only if the Riemann
Hypothesis is valid.

\section{Almost Counterexample to Riemann Hypothesis (Lehmer)}

An important landmark in the numerical study of the zeros of the zeta-function
is Lehmer's paper \cite{Le} in the year 1956. In it, he proves that the first
$10000$ zeros of the function have a real part exactly equal to $1/2$, so that
the Riemann Hypothesis is valid at least for $t\leq 9878.910$.

He establishes that, at this height, one out of ten Gram interval does not
satisfy Gram's law. The number of exceptions increases continuously. He also
finds out that in many occasions, in order to separate a zero, he must turn to
Euler-MacLaurin formula because Riemann-Siegel's one does not have enough
precission. This happens because  the zeros of the
zeta-function are very close. He studies specifically a particularly difficult
case: it is an area near $t=1114,89$, situated in the Gram interval $(g_{6707},
g_{6708})$, where the function has two extremely close zeros,
$${1\over2}+i\; 7005.0629\qquad {1\over2}+i\; 7005.1006.$$
We will not repeat the graphics that Lehmer made about the behaviour of the
function $Z(t)$ in a neighborhood of these points. Between these two zeros 
Hardy's function has the lowest relative maximum. This maximum is
only $0.0039675$ and it occurs at the point $t=7005.0819$. Looking at the
terms of $Z(t)$ in this point we see that a few of the first terms quickly
increasing are counteracted by conspiracy of lots of small terms which sum up,
thus the maximum turns out almost negative. A negative relative maximum would
imply, it can be proved, a counterexample to the Riemann Hypothesis. So we call
this situation an almost counterexample to the Riemann Hypothesis.

\vfil\eject

%%  MetaPost
%%  RADIOGRAFIA::05-zeta6990/7010::Zeta6990/7010.mp
%%

\end{large}
\begin{minipage}{330pt}
\begin{minipage}{154.5pt} 
\includegraphics[height=518.25pt]{figure37.ps}
\end{minipage}
\hfil
\begin{minipage}{154.5pt}
\includegraphics[height=518.25pt]{figure38.ps}
\end{minipage}
\centerline{\ }
\centerline{ {\scshape Figure} 17. Almost counterexample of Lehmer}
\end{minipage}

\begin{large}

\vfil\eject

In the figures included in this page we can see the case analyzed by Lehmer. We
notice that, from our point of view, it consists of two very close zeros, so
that, viewing it from a far distance, it seems to be a double zero. We see that
the lines seem to be continued more smoothly by the ones which are not actually
joining them. Instead, the lines turn abruptly. Every time this is the case,
that is, there are drastic changes in the direction of the lines, there is a
zero of the derivative hanging about. 

It would be easy to modify the lines
artificially so that the two lines containing the Gram points would join and
the other two thick lines would join each other too. So, slightly modifying the
path of the thin lines, we could generate two zeros outside the critical line,
which are symmetrical with respect to each other.

\end{large}

\begin{minipage}{325pt}
\bigskip
\includegraphics[width=325pt]{figure39.ps}
\centerline{\ }
\centerline{{\scshape Figure} 18. Detail of the last figure.}
\end{minipage}

\begin{large}

\section{Rosser Law}

Rosser, Yohe and Schoenfeld (1968) expand Lehmer's calculations and prove that
the first $3\,500\,000$ zeros are simple and situated on the critical line. These
authors find out a certain regularity in the failures of Gram's law.

They distinguish between {\sl good} and {\sl bad} Gram points. Good points are
those at which $\zeta(1/2 + ig_n)>0$ holds, and bad ones are the rest of them.
We know how to tell these points apart with the nacked eye in the X 
ray.
Bad Gram points are surrounded by a thin line.

They call {\bf Gram block} to a consecutive set of bad Gram points surrounded
by two good ones. For example, in the precedent figure points $g_{6707}$,
$g_{6708}$ and $g_{6709}$ form a Gram block.

Rosser's law claims that in a Gram block there are as many zeros as the number
of Gram's intervals.

\subsection{The function $S(t)$}
Let $N(T)$ be the number of zeros
$\rho = \beta + i\gamma$ with $0\leq \gamma \leq T$. An
approximation to $N(T)$ is $\pi^{-1}\theta(T) + 1$, so that
$$N(T)=\pi^{-1}\theta(T)+1+S(T).$$
The value  $\pi S(T)$  is also the variation of the argument of $\zeta(s)$
when $s$ goes from $+\infty+iT$ to $1/2+i T$. 

Von Mangoldt proved  that,  as Riemann says, $S(T)=\Orden(\log T)$. 
Later Littlewood proved 
$$\int_0^TS(t)\,dt=\Orden(\log T).$$

Selberg proves that $$S(t)=\Omega_\pm\bigl((\log t)^{1/3}(\log\log t)^{-7/3}\bigr).$$
Thus, there exist values of $t$ at which $S(t)$ is as high as we want it to be.

\subsection{First counterexample to Rosser's law.}
The first counterexample to
Rosser's law (see Figure $9$) is in the Gram block
$$(g_{13999525}, g_{13999527})$$  
in which there is no zero of the function. In the following interval
$J=(g_{13999527}, g_{13999528})$ there are three zeros which balance the total
count.

In the graphics we see how the function $S(T)$ takes a value greater than $2$
in a point which is situated between point $g_{13999527}$ and the first zero
contained in the interval $J$.

This is a general rule: Gram's law is satisfied as long as $\vert S\vert<1$ and
Rosser's law as long as $\vert S\vert <2$.

We have already said that there exist points on which $S$ takes values as high
as desired.

\vfil\eject

\end{large}

\begin{minipage}{330pt} 
\begin{minipage}{154.5pt}
\includegraphics[height=518.25pt]{figure40.ps}
\end{minipage}
\hfil
\begin{minipage}{154.5pt}
\includegraphics[height=518.25pt]{figure41.ps}
\end{minipage}
\centerline{\ }
\centerline{ {\scshape Figure} 19. First counterexample to Rosser's law.}
\end{minipage}

\begin{large}
\vfil\eject

A high value of $S(t)$ corresponds to a point in the critical
line, $1/2 + it$, such that the segment from $1/2 + it$ to
$+\infty + it$ meets a high number of our lines.

In the tables from \cite{B} and \cite{BLRW} we see that the extreme values of $S$ which are
known hardly surpass an absolute value of $2$.

In Brent's table is already quoted the first counterexample to
Rosser's law associated to the Gram point number $13999525$, where
$S(t)$ reaches a value of $-2.004138$. In fact, we see that there
is a point $1/2 + it_0$ imperceptibly above $g_{13999527}$, but before
the next zero of $\zeta(s)$, such that the segment from this point
to $1 + it_0$ meets firstly a thick line, then a thin line,
followed by another thick and another thin line. In the point
$+\infty + it_0$ we start from a value equal to $1$, enter the
fourth quadrant, cross a thin line and thus we reach the third
quadrant, cross a thick line entering the second quadrant, cross
another thin line and find ourselves in the first quadrant,  and
finally another thick line and so we end up in the fourth
quadrant. Consequently, the argument has changed in a quantity
between $2\pi$ and $5\pi/2$, that, as we can see, agrees with the
value given by Brent.

\vfil
\end{large}
\begin{minipage}{321pt}
\bigskip
\includegraphics[width=312pt]{figure42.ps}
\vfil
\centerline{{\scshape Figure} 20. Detail of the last figure}
\end{minipage}
\goodbreak

\begin{minipage}{330pt} 
\begin{minipage}{154.5pt}
\includegraphics[height=518.25pt]{figure43.ps}
\end{minipage}
\hfil
\begin{minipage}{154.5pt}
\includegraphics[height=518.25pt]{figure44.ps}
\end{minipage}
\centerline{\ }
\centerline{ {\scshape Figure} 21. Another counterexample to  Rosser's law.}
\end{minipage}

\vfil\eject

\begin{minipage}{321pt}
\bigskip
\includegraphics[width=312pt]{figure45.ps}
\vfil
\centerline{{ \scshape Figure} 22. Detail of the last figure.}
\end{minipage}
\begin{large}
\vskip0.2cm
In the example presented in this page, $S$ reaches a value $2.0506$ and so the
existence of the point $1/2 + it_0$ is a little clearer. In this case, $S$ is
positive, while in the precedent it was negative.

The previously quoted result of Selberg assures us that there are points at
which $S(t)$ is as high as desired. Thus we can assume that, in a higher level,
we will see coils which are analogous to Figures 22 and 20, but in
which an arbitary number of lines gets involved.

Watching the preceding figures one may wonder if the thin lines 
do not cross the line $\sigma=0$. This is not true, by a theorem of Bohr (see 
\cite{Ti} p.~300) 
the function $\zeta(s)$ takes every value $\ne 0$ infinitely often in the 
halfplane 
$\Re(s)>1$. In fact Van de Lune \cite{Lu} has shown that 
$\sigma_0=\sup\{\sigma\in\R: \Re\zeta(\sigma+it)<0$ for some $t\in\R\}$ is given by 
the unique solution of the equation $\sum_p \arcsin(p^{-\sigma})=\pi/2$, $\sigma>1$. 
Brent and Van de Lune have computed 
$\sigma_0=1.1923473371861932\dots$ with more than 400 decimal digits. 
%Also van de Lune proved that the thin lines do not reach the line $\Re(s)=\sigma_0$.
\par
We finish with three figures of the zeta-function near a thousand millions, 
to show a ramdomly chosen area.
\goodbreak

\end{large}

\begin{minipage}{335.7426pt}
\hskip-24.2574pt\includegraphics[width=360pt]{figure46.ps} 
\end{minipage}
\centerline{\ }
\vfil
\centerline{{ \scshape Figure} 23. $\zeta(s)$ near $t=1\,000\,000\,000$.}
\vfil\eject

\begin{minipage}{335.7426pt}
\hskip-24.2574pt\includegraphics[width=360pt]{figure47.ps} 
\end{minipage}
\centerline{\ }
\vfil
\centerline{{\scshape Figure 24} .  $\zeta(s)$ near $t=1\,000\,000\,000$.}
\vfil\eject

\begin{minipage}{335.7426pt}
\hskip-24.2574pt\includegraphics[width=360pt]{figure48.ps} 
\end{minipage}
\centerline{\ }
\vfil
\centerline{{\scshape Figure} 25.  $\zeta(s)$ near $t=1\,000\,000\,000$.}
\vfil\eject

\begin{minipage}{321pt}
\bigskip
\includegraphics[width=312pt]{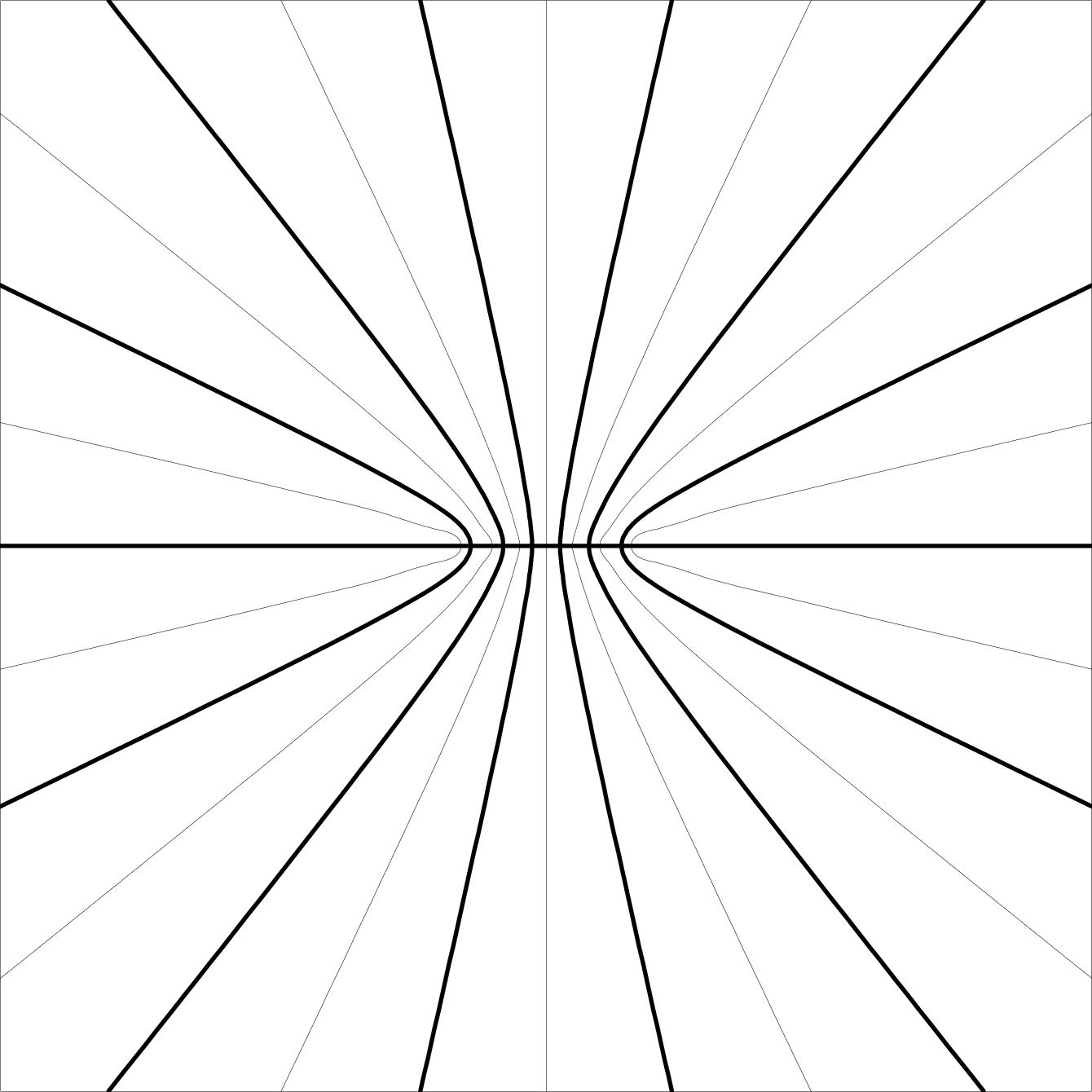}
\vfil
\centerline{{\scshape Figure} 26. Hermite polynomial $H_{7}(z)$}
\end{minipage}
\medskip

\begin{large}
    
    To show how boring life can be outside Number Theory, we include the graphics
of some functions.

The first one is Hermite's polynomial $H_7(z)$, that is,
$$   128 z^7- 1344 z^5+ 3360 z^3-1680 z $$
It has a degree equal to seven and all its roots are real. Here we represent it
in the rectangle $(-17, 17)^2$, which is enough to get a clear idea of how the
graphic is.

We can see the seven zeros of the function and the six ceros of the derivative.
\bigskip

The graphics of all Hermite's polynomials are analogous. It also looks 
as the
graphics of other orthogonal polynomial families. But we must point out that a
general polynomial can have a very complicated graphics. The regularity in this
case is due to the fact that it is a very particular polynomial.

In this page we have the X rays corresponding to the Bessel
function $J_7(z)$ and the Airy function $Ai(z)$ defined by

\end{large}
\begin{minipage}{335.7426pt}
\bigskip
\includegraphics[width=335pt]{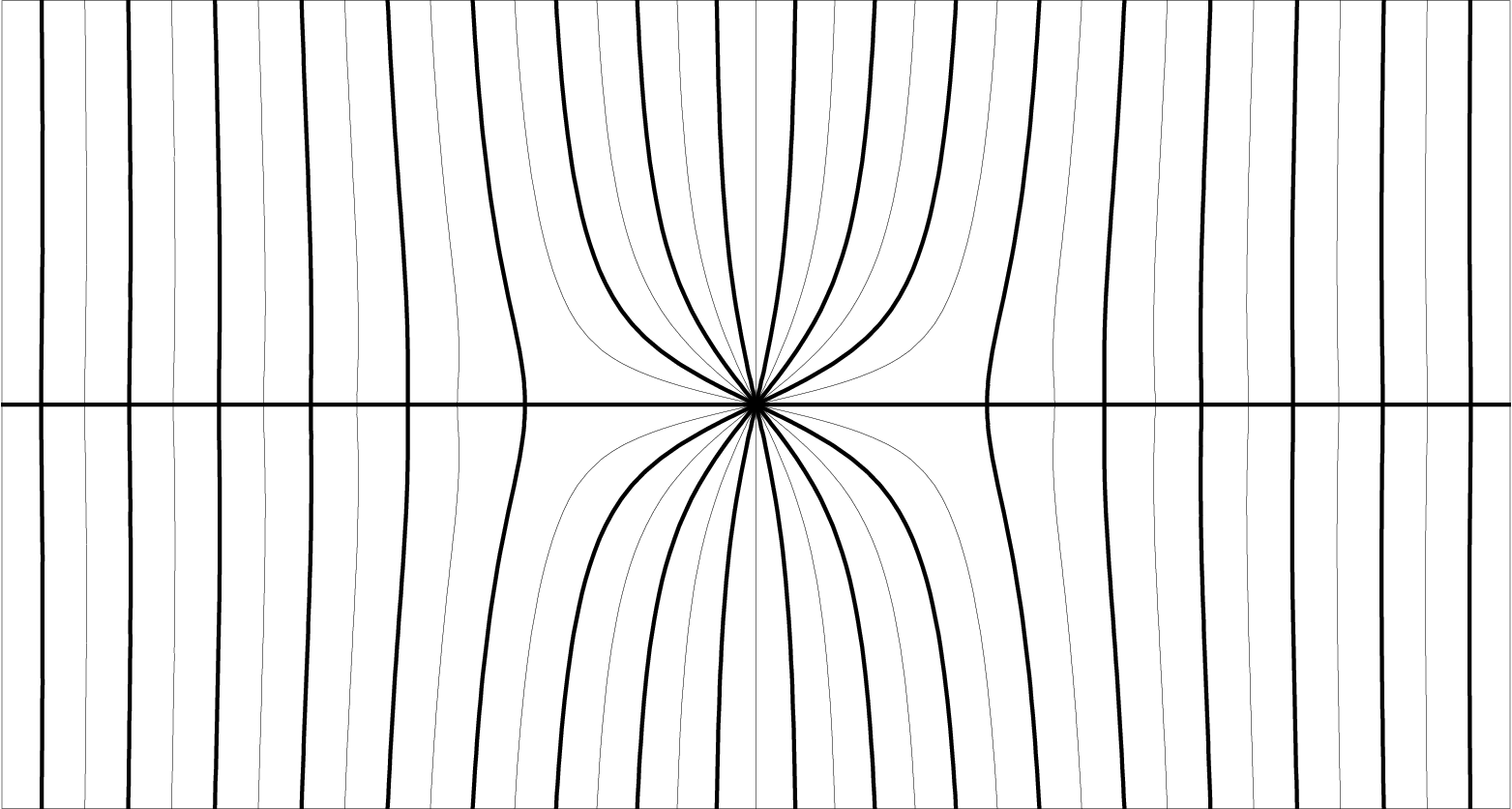}
\vfil
\centerline{{\scshape Figure} 27. Bessel function $J_{7}(z)$}
\bigskip
\includegraphics[width=335pt]{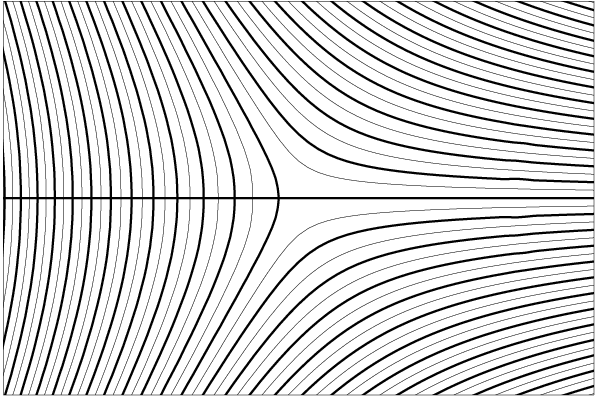}
\vfil
\centerline{{\scshape Figure} 28.  Airy function $Ai(z)$}
\end{minipage}

\begin{large}

$$\mskip-20mu J_{7}(z)=\Bigl({z\over2}\Bigr)^7\sum_{k=0}^\infty {(-1)^k 
(z/2)^k\over k! (k+7)!}, Ai(z)={3^{-2/3}\over 
\pi}\sum_{k=0}^\infty{\Gamma\bigl((k+1)/3\bigr)\sen{2\pi\over3}(k+1)\over 
k!}(3^{1/3}z)^k\mskip-70mu$$

\bigskip

The Bessel function is represented on the rectangle
 $(-28, 28)\times (-20, 20)$ and the Airy function on  $(-15, 15)\times (-10,
10)$. 

The Bessel function has a zero of order $7$ in the origen. Its other
zeros are real and, appart from the obvious zeros of the derivative, which are
real, the derivative does not vanish.

The Airy function is surprising because of the likeness of its X ray 
with that of the
Gamma function.

\end{large}
\begin{minipage}{321pt}
\bigskip
\includegraphics[width=312pt]{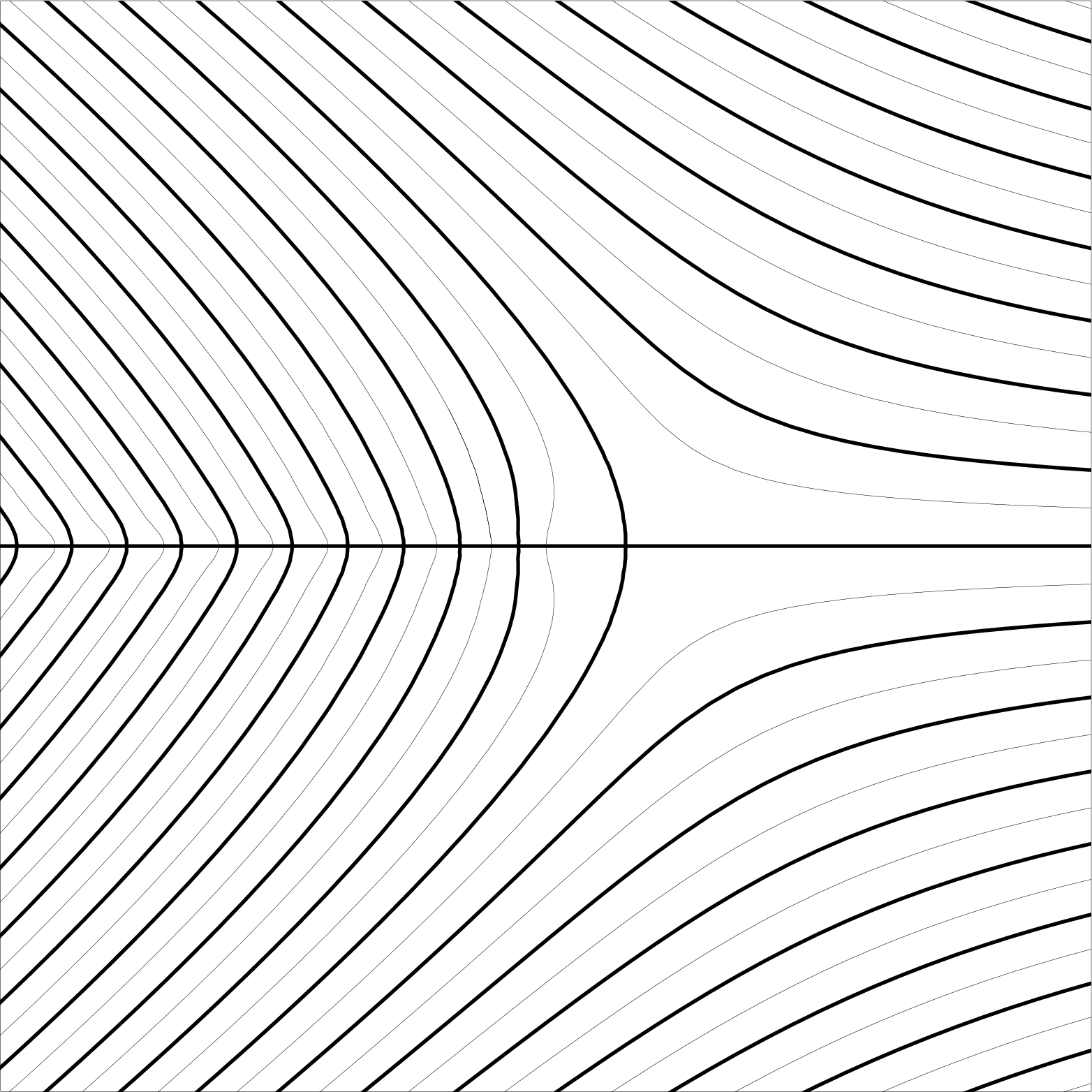}
\vfil
\centerline{{\scshape Figure} 29. Function $\Gamma(s)$}
\end{minipage}
\medskip

\begin{large}

The graphic of the function $\Gamma(s)$ shows that its derivative vanishes
only at the obvious zeros. The graphics of a function and its inverse do always
coincide. The figure shows the rectangle $(-10, 10)^2$.

\end{large}

\vspace{3cm}

\noindent  This paper was presented on a lecture on Analytic Number Theory on 
   May 29, 2002

\vspace{4cm}

\end{document}